\documentclass{article}%
\usepackage{amsfonts}
\usepackage{amssymb}
\usepackage{amsmath}
\usepackage{graphicx}%
\providecommand{\U}[1]{\protect\rule{.1in}{.1in}}

\begin{document}

\title{\textbf{The Logical Theory of Canonical Maps}: \\The Elements \& Distinctions Analysis of the Morphisms, Duality, Canonicity,
and Universal Constructions in $Sets$}
\author{David Ellerman orcid.org/0000-0002-5718-618X\\Independent Researcher\\Email: david@ellerman.org}
\date{}
\maketitle

\begin{abstract}
\noindent Category theory gives a mathematical characterization of naturality
but not of canonicity. The purpose of this paper is to develop the logical
theory of canonical maps based on the broader demonstration that the dual
notions of elements \& distinctions are the basic analytical concepts needed
to unpack and analyze morphisms, duality, canonicity, and universal
constructions in $Sets$, the category of sets and functions. The analysis
extends directly to other $Sets$-based concrete categories (groups, rings,
vector spaces, etc.). Elements and distinctions are the building blocks of the
two dual logics, the Boolean logic of subsets and the logic of partitions. The
partial orders (inclusion and refinement) in the lattices for the dual logics
define the canonical morphisms (where `canonical' is always relative to the
given data, not an absolute property of a morphism). The thesis is that the
maps that are canonical in $Sets$ are the ones that are defined (given the
data of the situation) by these two logical partial orders and by the
compositions of those maps.


AMS: 03, 18

\end{abstract}

\section{Elements \& Distinctions Analysis}

\subsection{Introduction}

Category theory gives a mathematical characterization of naturality
\cite{eilen-macl:nat-equi} but not of canonicity, the canonical nature of
certain morphisms. The purpose of this paper is to present the \textit{logical
theory of canonical maps} that provides such a characterization. That logical
theory of canonical maps is one of the main results in the broader analysis
showing that the dual notions of "elements \& distinctions" (or "its \& dits")
are the basic analytical concepts needed to unpack and analyze morphisms,
duality, canonicity (or canonicalness), and universal constructions in $Sets$,
the category of sets and functions. The analysis extends directly to other
$Sets$-based concrete categories (groups, rings, vector spaces, etc.) where
the objects are sets with a certain type of structure and the morphisms are
set functions that preserve or reflect that structure. Then the elements \&
distinctions-based definitions can be abstracted in purely arrow-theoretic way
for abstract category theory.

One way to approach the concepts of "elements" (or "its") and "distinctions"
(or "dits") is to start with the category-theoretic duality between subsets
and quotient sets (= partitions = equivalence relations): "The dual notion
(obtained by reversing the arrows) of `part' [subobject] is the notion of
\textit{partition}."\ \cite[p. 85]{law:sfm}. That motivates the two dual forms
of mathematical logic: the Boolean logic of subsets and the logic of
partitions (\cite{ell:lop}; \cite{ell:intropartitions}). If partitions are
dual to subsets, then what is the dual concept that corresponds to the notion
of \textit{elements of a subset}? The notion dual to the elements of a subset
is the notion of the \textit{distinctions of a partition} (pairs of elements
in distinct blocks of the partition).

\subsection{Set functions transmit elements and reflect distinctions}

The duality between elements ("its") of a subset and distinctions ("dits") of
a partition already appears in the very notion of a function between sets in
the category $Sets$. The concepts of elements and distinctions provide the
natural notions to specify the binary relations, i.e., subsets $R\subseteq
X\times Y$, that define functions $f:X\rightarrow Y$.

A binary relation $R\subseteq X\times Y$ \textit{transmits elements} if for
each element $x\in X$, there is an ordered pair $\left(  x,y\right)  \in R$
for some $y\in Y$.

A binary relation $R\subseteq X\times Y$ \textit{reflects elements} if for
each element $y\in Y$, there is an ordered pair $\left(  x,y\right)  \in R$
for some $x\in X$.

A binary relation $R\subseteq X\times Y$ \textit{transmits distinctions} if
for any pairs $\left(  x,y\right)  $ and $\left(  x^{\prime},y^{\prime
}\right)  $ in $R$, if $x\not =x^{\prime}$, then $y\not =y^{\prime}$.

A binary relation $R\subseteq X\times Y$ \textit{reflects distinctions} if for
any pairs $\left(  x,y\right)  $ and $\left(  x^{\prime},y^{\prime}\right)  $
in $R$, if $y\not =y^{\prime}$, then $x\not =x^{\prime}$.

The dual role of elements and distinctions can be seen if we translate the
usual characterization of the binary relations that define functions into the
elements-and-distinctions language. In the usual treatment, a binary relation
$R\subseteq X\times Y$ defines a \textit{function} $X\rightarrow Y$ if it is
defined everywhere on $X$ and is single-valued. But "being defined everywhere"
is the same as transmitting (or "preserving") elements, and being
single-valued is the same as reflecting distinctions so the more natural
definition is:

\begin{center}
a binary relation $R$ is a \textit{function} if it transmits elements and
reflects distinctions.
\end{center}

What about the other two special types of relations, i.e., those which
transmit (or preserve) distinctions or reflect elements? The two important
special types of functions are the injections and surjections, and they are
defined by the other two notions:

\begin{center}
a function is \textit{injective} if it transmits distinctions, and

a function is \textit{surjective} if it reflects elements.
\end{center}

Given a set function $f:X\rightarrow Y$ with \textit{domain} $X$ and
\textit{codomain} $Y$, a subset of the codomain $Y$ is determined as the
\textit{image} $f(X)\subseteq Y$, and a partition on the domain $X$ is
determined as the \textit{coimage} or \textit{inverse-image} $\{f^{-1}%
(y)\}_{y\in f\left(  X\right)  }$. It might also be noted that the empty set
of ordered pairs $\emptyset\times Y$ satisfies the definition of a function
$\emptyset\rightarrow Y$ whose image is the empty subset $\emptyset$ of $T$
and coimage is the empty partition $\emptyset$ on $\emptyset$.

\subsection{The logical theory of \textit{canonical} maps based on the its \&
dits analysis}

Jean-Pierre Marquis \cite{marquis:canonical} has raised the question of
characterizing canonical maps in mathematics in general and category theory in
particular. Category theory gives a mathematical notion of "naturality" but
not of canonicalness or canonicity. Marquis gives the intuitive idea (maps
defined "without any arbitrary decision"), a number of examples (most of which
we will analyze in $Sets$), and a set of criteria stated in terms of limits
(and thus dually for colimits).

\begin{quotation}
We are now in a position to circumscribe more precisely what we want to
include in the notion of canonical morphisms or maps.

1. Morphisms that are part of the data of a limit are canonical morphisms; for
instance, the projection morphisms that are part of the notion of a product;

2. The unique morphism from a cone to a limit determined by a universal
property is a canonical morphism: and

3. In particular, the unique isomorphism that arise between two candidates for
a limit is a canonical morphism. \cite[p. 101]{marquis:canonical}
\end{quotation}

The elements \& distinctions (or its \& dits) analysis provides a
\textit{mathematical} characterization of "canonical maps" in $Sets$ (and thus
in $Sets$-based concrete categories) that satisfies the Marquis criteria. The
characterization of canonicity is always \textit{relative to the given data};
canonicity is not an `absolute' property of a morphism. For instance, given a
cone $f:Z\rightarrow X$ and $g:Z\rightarrow Y$, the canonical map
$Z\rightarrow X\times Y$ is only canonical relative to the data $f$ and $g$.

The treatment of canonicity is part of a broader elements \& distinctions
analysis of the morphisms, duality, and universal constructions (limits and
colimits) in the basic `ur-category' $Sets$ of sets and functions and thus in
$Sets$-based categories--which is abstracted in category theory as a whole. At
that point, the elements \& distinctions analysis connects to the broader
philosophical literature on structuralism since category theory is essentially
the natural codification of that philosophy of mathematics
(\cite{landry:structuralism}; \cite{landrymarguis:catincontext};
\cite{landry:catsworkingphil}).

The logical theory of canonicity is that the canonical maps and the unique
canonical factor morphisms in the universal mapping properties (UMPs) in
$Sets$ are always constructed in the two ways that maps are constructed
\textit{from the partial orders} in the two basic logics, the logic of subsets
and the logic of partitions. In the powerset Boolean algebra of subsets
$\wp\left(  U\right)  $ of $U$, the partial order is the inclusion relation
$S\subseteq T$ for $S,T\subseteq U$, which induces the canonical injection
$S\rightarrowtail T$. That is the way canonical injective maps are defined
from the partial order of inclusion on subsets.

In the dual algebra of partitions $\Pi\left(  U\right)  $ on $U$ (i.e., the
lattice of partitions on $U$ enriched with the implication operation on
partitions\footnote{In \cite{ell:intropartitions}, the partition algebra was
defined as the partition lattice enriched with the implication and nand
operations on partitions. But for purposes of comparisons with Boolean or
Heyting algebras, it suffices to consider only the implication in addition to
the join and meet. In any case, this does not affect the analysis here where
the lattice structure suffices.}), the partial order is the refinement
relation between partitions and it induces a canonical map using refinement. A
\textit{partition} $\pi=\left\{  B,B^{\prime},...\right\}  $ on a set $U$ is a
set of non-empty subsets of $U$ (called \textit{blocks}, $B,B^{\prime},...$)
that are mutually exclusive (i.e., disjoint) and jointly exhaustive (i.e.,
whose union is $U$). It might be noticed that the empty set $\emptyset$, which
has no nonempty subsets, is the empty partition on $U=\emptyset$%
.\footnote{Thanks to Paul Blain Levy and Alex Simpson for emphasizing to me
the role of empty partition for the consistent development of the whole
theory, e.g., as the inverse-image partition on the domain of the empty
function $\emptyset\rightarrow Y$.} One could also define a partition on $U$
as the set of inverse-images $f^{-1}\left(  y\right)  \subseteq U$ for $y\in
f\left(  U\right)  $ for any \textit{function} $f:U\rightarrow Y$ with domain
$U$. The empty partition on $U=\emptyset$ is then the inverse-image partition
of the empty function $\emptyset\rightarrow Y$.

Given another partition $\sigma=\left\{  C,C^{\prime},...\right\}  $ on $U$, a
partition $\pi$ is said to \textit{refine }$\sigma$ (or $\sigma$ is
\textit{refined by} $\pi$), written $\sigma\precsim\pi$, if for every block
$B\in\pi$, there is a block $C\in\sigma$ (necessarily unique) such that
$B\subseteq C$. If we denote the set of distinctions or dits of a partition
(ordered pairs of elements in different blocks) by $\operatorname{dit}\left(
\pi\right)  $, the ditset of $\pi$, then just as the partial order in
$\wp\left(  U\right)  $ is the inclusion of elements, so the refinement
partial order on $\Pi\left(  U\right)  $ is the inclusion of distinctions,
i.e., $\sigma\precsim\pi$ iff (if and only if) $\operatorname{dit}\left(
\sigma\right)  \subseteq\operatorname{dit}\left(  \pi\right)  $. And just as
the inclusion ordering on subsets induces a canonical injection between
subsets, so the refinement ordering $\sigma\precsim\pi$ on partitions induces
a canonical surjection between partitions, namely $\pi\rightarrow\sigma$ where
$B\in\pi$ is taken to the unique $C\in\sigma$ where $B\subseteq C$. If the
blocks of $\pi=\left\{  B_{x}\right\}  _{x\in X}$ are indexed by a set $X$ and
the blocks of $\sigma=\left\{  C_{y}\right\}  _{y\in Y}$ are indexed by a set
$Y$, then the refinement $\sigma\precsim\pi$ induces a canonical surjection
$X\twoheadrightarrow Y$. That is the way canonical surjective maps are defined
from the partial order of refinement on partitions.

These canonical injections and surjections are \textit{built into} the partial
orders of the lattice (or algebraic) structure of the two dual logics of
subsets and partitions; they logically define the `atomic' \textit{canonical}
maps in $Sets$, and other canonical maps in $Sets$ arise out of their
compositions. Note that the canonical injections are in the `upward' direction
of the partial order (more elements) in the lattice of subsets, while the the
canonical surjections are in the opposite downward direction to the `upward'
(more dits) direction of the partial order in the lattice of partitions.

This \textit{logical theory of canonical maps} is that all "canonical" maps in
$Sets$ arise from the given data in these two ways or by compositions of
them--which then extends to $Sets$-based concrete categories. "Canonical"
always means \textit{relative} to the given data. The given data only plays
the role of defining the sets with the inclusion relations between them or the
partitions with the refinement relation between them. That is the role of the
given data. Then the canonical injections and canonical surjections are
defined by those inclusions or refinements. Marquis' informal definition of
canonical maps as maps defined "without any arbitrary decision" then means
that once the initial data is encoded as inclusion or refinement relations in
the respective logical lattices, then the logical structure suffices to define
the canonical injective or surjective maps. And the thesis is that all
canonical maps in $Sets$ are those canonical injections, canonical
surjections, or their compositions.

The thesis cannot be proven since "canonical" is an intuitive notion. But we
will show that all the canonical maps and unique factor maps in the universal
constructions (limits, colimits, and the exponential or Currying adjuction) in
$Sets$ arise in this way from the partial orders of the dual lattices (or
algebras) of subsets and partitions--which thus satisfies the Marquis
criteria. This logical basis for this theory of canonical maps accounts for
the name "logical." 

\subsection{Initial \& terminal objects in $Sets$}

The top of the powerset Boolean algebra $\wp\left(  U\right)  $ is $U$, where
each subset $S\subseteq U$ induces the canonical injection $S\rightarrow U$.
The bottom of the Boolean algebra, the null set $\emptyset$, is included in
any set, e.g., $\emptyset\subseteq U$, so the induced morphism $\emptyset
\rightarrow U$ is the canonical map that makes $\emptyset$ the initial object
in $Sets$ (taking $U$ as any set).

The top of the partition algebra $\Pi\left(  U\right)  $ is the
\textit{discrete partition} $\mathbf{1}_{U}=\left\{  \left\{  u\right\}
\right\}  _{u\in U}$ of all singletons. Since every partition $\pi$ is refined
by $\mathbf{1}_{U}$, i.e., $\pi\precsim\mathbf{1}_{U}$, there is the canonical
surjection $\mathbf{1}_{U}\cong U\rightarrow\pi$ that takes the singleton
$\left\{  u\right\}  $ or just $u$ (since blocks of $\mathbf{1}_{U}$ are in
one-to-one correspondence with the elements of $U$) to the unique block $B$
such that $u\in B$. The bottom of the partition algebra (or lattice) is the
\textit{indiscrete partition} (nicknamed the "Blob") $\mathbf{0}_{U}=\left\{
U\right\}  $ with only one block $U$ that identifies all the points in $U$ so
$\mathbf{0}_{U}$ is isomorphic to the one-element set $1$. And $\mathbf{0}%
_{U}$ is refined by all partitions, e.g., $\mathbf{0}_{U}\precsim
\mathbf{1}_{U}$. That refinement relation induces the unique map from the
blocks of $\mathbf{1}_{U}$ (i.e., the elements of $U$) to the blocks or rather
\textquotedblleft the block\textquotedblright\ of $\mathbf{0}_{U}\cong1$,
i.e., induces the canonical map $U\twoheadrightarrow1$ that makes the
one-element set $1$ into the terminal object in $Sets$ (taking $U$ as any set).

Thus the maps induced by the bottom of the inclusion/refinement relations in
the two logical partial orders give the canonical maps for the initial and
terminal objects in $Sets$.

\begin{center}%
\begin{tabular}
[c]{|c||c|c|}\hline
Dualities & Subset logic & Partition logic\\\hline\hline
`Elements' & Elements $u$ of $S$ & Dits $\left(  u,u^{\prime}\right)  $ of
$\pi$\\\hline
Partial order & Inclusion $S\subseteq T$ & $\sigma\precsim\pi$:
$\operatorname*{dit}\left(  \sigma\right)  \subseteq\operatorname*{dit}\left(
\pi\right)  $\\\hline
Canonical map & $S\rightarrowtail T$ & $\pi\twoheadrightarrow\sigma$\\\hline
Top of partial order & $U$ all elements & $\mathbf{1}_{U}$,
$\operatorname*{dit}(\mathbf{1}_{U}\mathbf{)}=U^{2}-\Delta$, all dits\\\hline
Bottom of partial orders & $\emptyset$ no elements & $\mathbf{0}_{U}$,
$\operatorname*{dit}(\mathbf{0}_{U}\mathbf{)=\emptyset}$, no dits\\\hline
Extremal objects $Sets$ & $\emptyset\subseteq U$, $\emptyset\rightarrowtail U$
& $1\cong\mathbf{0}_{U}\precsim\mathbf{1}_{U}$, $U\twoheadrightarrow1$\\\hline
\end{tabular}

Table 1: Elements and distinctions in the dual logics
\end{center}

There are different ways to characterize objects with universal mapping
properties (like the initial and terminal objects) using `higher order'
machinery in category theory, e.g., using units or counits of adjunctions or
using representable functors. For instance, consider the covariant functor
$I:\mathcal{C}\rightarrow Sets$ that takes any object in a category
$\mathcal{C}$ to a singleton set in $Sets$. If there is a natural isomorphism
$Hom_{\mathcal{C}}(0,-)\cong I\left(  -\right)  $ then the object
$0\in\mathcal{C}$ \textit{represents} $I\left(  -\right)  $ and $0$ is an
initial object in $C$. If there is another object $0^{\prime}$ such that
$Hom_{\mathcal{C}}(0^{\prime},-)\cong I\left(  -\right)  $, then
$Hom_{\mathcal{C}}(0,-)\cong Hom_{\mathcal{C}}(0^{\prime},-)$ and, by the
Yoneda Lemma, there is a canonical isomorphism $0\cong0^{\prime}$
(\cite{biss:functor}; \cite{mazur:equal}). But this important characterization
of the initial object $\emptyset$ in $Sets$ does not offer an explanation of
why $\emptyset$ has that universal property in the first place. Our claim is
that the underlying fact that $\emptyset\subseteq U$ in $\wp\left(  U\right)
$ for any set $U$ accounts for it being the initial object in $Sets$, and
dually, $1\cong\mathbf{0}_{U}\precsim\mathbf{1}_{U}$ in $\Pi(U)$ for any set
$U$ accounts for $1$ being the terminal object in $Sets$. And $\emptyset$
being the no-elements subset and $\mathbf{0}_{U}$ being the no-distinctions
partition accounts for them being the bottoms of the dual logical lattices
$\wp\left(  U\right)  $ and $\Pi(U)$ where the respective partial orders are
inclusions of elements and inclusions of distinctions.

\subsection{The epi-mono factorization in $Sets$}

Another simple application of the elements \& distinctions analysis is the
construction of the canonical surjection and canonical injection in the
epi-mono factorization of any set function: $f:X\rightarrow Y$. The data in
the function provide the \textit{coimage} (or inverse-image) partition
$f^{-1}=\left\{  f^{-1}\left(  y\right)  :y\in f\left(  X\right)  \right\}  $
on $X$ [where the blocks of $f^{-1}$ are indexed by the $y\in f\left(
X\right)  $] and the \textit{image} subset $f\left(  X\right)  $. Since
$f^{-1}$ is refined by the discrete partition on $X$, $f^{-1}\precsim
\mathbf{1}_{X}$, the induced surjection is the canonical map
$X\twoheadrightarrow f\left(  X\right)  $. Note that in this case, the initial
data $f$ only defined one partition $f^{-1}$ on $X$. The refinement relation
used is that the discrete partition $\mathbf{1}_{X}$ on $X$ refines all
partitions on $X$. Similarly the initial data $f$ only defines one subset
$f\left(  X\right)  $ of $Y$ but all such subsets are included in $Y$ so that
inclusion $f\left(  X\right)  \subseteq Y$ induces the injection $f\left(
X\right)  \rightarrowtail Y$ and the epi-mono factorization of $f$ is:

\begin{center}
$f:X\rightarrow Y=X\twoheadrightarrow f\left(  X\right)  \rightarrowtail Y$.
\end{center}

\subsection{Abstracting to arrow-theoretic definitions}

One of our themes is that the concepts of elements and distinctions unpack and
analyze the category theoretic concepts in the basic `ur-category' $Sets$, and
they are abstracted into purely arrow-theoretic definitions in abstract
category theory. For instance, the elements \& distinctions definitions of
injections and surjections yield "arrow-theoretic" characterizations which can
then be applied in any category to provide the usual category-theoretic dual
definitions of monomorphisms (injections for set functions) and epimorphisms
(surjections for set functions).

Two set functions $f,g:X\rightrightarrows Y$ are different, i.e., $f\not =g$,
if there is an element $x$ of $X$ such that their values $f\left(  x\right)  $
and $g\left(  x\right)  $ are a distinction of $Y$, i.e., $f\left(  x\right)
\not =g\left(  x\right)  $. Hence if $f$ and $g$ are followed by a function
$h:Y\rightarrow Z$, then the compositions $hf,hg:X\rightarrow Y\rightarrow Z$
must be different if $h$ \textit{preserves} \textit{distinctions} (so that the
distinction $f(x)\not =g\left(  x\right)  $ is preserved as $hf\left(
x\right)  \not =hg\left(  x\right)  $), i.e., if $h$ is injective. Thus in the
category of sets, $h$ being injective is characterized by: for any
$f,g:X\rightrightarrows Y$, "$f\neq g$ implies $hf\neq hg$" or equivalently,
"$hf=hg$ implies $f=g$" which is the general category-theoretic definition of
a \textit{monomorphism} or \textit{mono}.

In a similar manner, if we have functions $f,g:X\rightrightarrows Y$ where
$f\neq g$, i.e., where there is an element $x$ of $X$ such that their values
$f\left(  x\right)  $ and $g\left(  x\right)  $ are a distinction of $Y$, then
suppose the functions are preceded by a function $h:W\rightarrow X$. Then the
compositions $fh,gh:W\rightarrow X\rightarrow Y$ must be different if $h$
\textit{reflects elements} (so that the element $x$ where $f$ and $g$ differ
is sure to be in the image of $h$), i.e., if $h$ is surjective. Thus in the
category of sets, $h$ being surjective is characterized by: for any
$f,g:X\rightrightarrows Y$, "$f\neq g$ implies $fh\neq fg$" or "$fh=gh$
implies $f=g$" which is the general category-theoretic definition of an
\textit{epimorphism} or \textit{epi}.

Hence the dual interplay of the notions of elements \& distinctions can be
seen as yielding the arrow-theoretic characterizations of injections and
surjections which are lifted into the general categorical dual definitions of
monomorphisms and epimorphisms.

\subsection{Duality interchanges elements \& distinctions}

In the duality of plane projective geometry, every proof of a theorem
involving points and lines yields another proof of the theorem with points and
lines interchanges. Similarly, any arrow-theoretic proof of a result in
category theory yields a proof of a result in the opposite category with the
arrows rseversed. In $Sets$, what is interchanged (like points and lines) to
reverse the arrows? The reverse-the-arrows duality of category theory is the
abstraction from the reversing of the roles of elements \& distinctions (or
its \& dits) in dualizing $Sets$ to $Sets^{op}$. That is, a concrete morphism
in $Sets^{op}$ is a binary relation, which might be called a
\textit{cofunction}, that preserves distinctions and reflects
elements--instead of preserving elements and reflecting distinctions. %

\begin{center}
\includegraphics[
height=0.8389in,
width=4.1857in
]%
{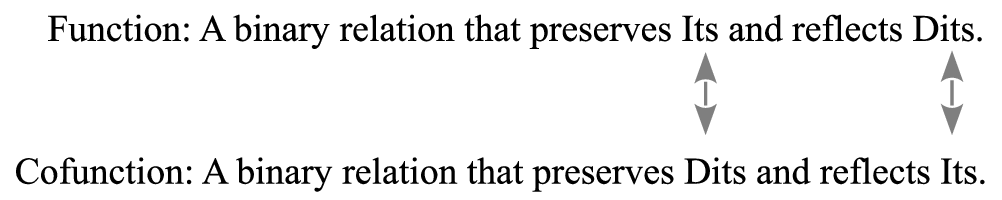}%
\end{center}

\begin{center}
Figure 1:\ Interchange Its \& Dits dualizes between functions and cofunctions
\end{center}

Equivalently, when we reverse the direction of a binary relation defining a
function, we just interchanged "reflects" and "preserves" (or "transmits").
Thus with every binary relation $f\subseteq X\times Y$ that is a function
$f:X\rightarrow Y$, there is a binary relation $f^{op}\subseteq Y\times X$
that is a cofunction $f^{op}:Y\rightarrow X$ in the opposite direction.%

\begin{center}
\includegraphics[
height=0.8354in,
width=4.1701in
]%
{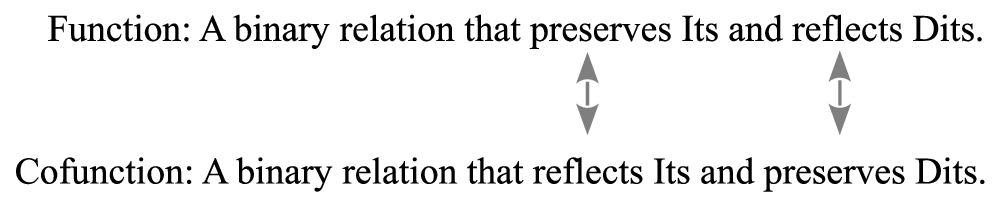}%
\end{center}

\begin{center}
Figure 2:\ Interchange preserves and reflects dualizes between functions and cofunctions
\end{center}

For the universal constructions in $Sets$, the interchange in the roles of
elements and distinctions interchanges each construction and its dual:
products and coproducts, equalizers and coequalizers, and in general limits
and colimits. That is then abstracted to make the reverse-the-arrows duality
in abstract category theory.

This begins to illustrate our theme that the language of elements \&
distinctions is the conceptual language in which the category of sets and
functions is written, and abstract category theory gives the abstract-arrows
version of those definitions. Hence we turn to universal constructions for
further analysis.

\section{The Elements \& Distinctions Analysis of Products and Coproducts}

\subsection{The coproduct in $Sets$}

Given two sets $X$ and $Y$ in $Sets$, the idea of the \textit{coproduct} is to
create the set with the maximum number of \textit{elements} starting with $X$
and $Y$. Since $X$ and $Y$ may overlap, we must make two copies of the
elements in the intersection. Hence the relevant operation is not the union of
sets $X\cup Y$ but the disjoint union $X\sqcup Y$. To take the disjoint union
of a set $X$ with itself, a copy $X^{\ast}=\left\{  x^{\ast}:x\in X\right\}  $
of $X$ is made so that $X\sqcup X$ can be constructed as $X\cup X^{\ast}$. In
a similar manner, if $X$ and $Y$ overlap, then $X\sqcup Y=X\cup Y^{\ast}$.
Then the inclusions $X,Y\subseteq X\sqcup Y$, give the canonical injections
$i_{X}:X\rightarrow X\sqcup Y$ and $i_{Y}:Y\rightarrow X\sqcup Y$.

The universal mapping property for the coproduct in $Sets$ is that given any
other `cocone' of maps $f:X\rightarrow Z$ and $g:Y\rightarrow Z$, there is a
unique map $f\sqcup g:X\sqcup Y\rightarrow Z$ such that the triangles commute
in the following diagram.

\begin{center}
$%
\begin{array}
[c]{ccccc}%
X & \overset{i_{X}}{\longrightarrow} & X\sqcup Y & \overset{i_{Y}%
}{\longleftarrow} & Y\\
& \searrow^{f} & ^{\exists!}\downarrow^{f\sqcup g} & ^{g}\swarrow & \\
&  & Z &  &
\end{array}
$

Coproduct diagram
\end{center}

From the data $f:X\rightarrow Z$ and $g:Y\rightarrow Z$, we need to
canonically construct the unique factor map $f\sqcup g:X\sqcup Y\rightarrow
Z$. The map $f:X\rightarrow Z$ defines the image $f\left(  X\right)  \subseteq
Z$ and $g:Y\rightarrow Z$ defines the image $g\left(  Y\right)  \subseteq Z$
and the subset lattice join in $\wp\left(  Z\right)  $ gives $f\left(
X\right)  \cup g\left(  Y\right)  \subseteq Z$ so there is a canonical
injection $f\left(  X\right)  \cup g\left(  Y\right)  \rightarrowtail Z$.
Since each $x\in X\sqcup Y$ is associated with $f\left(  x\right)  \in
f\left(  X\right)  \cup g\left(  Y\right)  $ and each $y\in X\sqcup Y$ is
associated with $g\left(  y\right)  \in f\left(  X\right)  \cup g\left(
Y\right)  $, that epi completes the factor map $f\sqcup g:X\sqcup
Y\twoheadrightarrow f\left(  X\right)  \cup g\left(  Y\right)  \rightarrowtail
Z$.

\subsection{The product in $Sets$}

Given two sets $X$ and $Y$ in $Sets$, the idea of the \textit{product} is to
create the set with the maximum number of \textit{distinctions} starting with
$X$ and $Y$. The product in $Sets$ is usually constructed as the set of
ordered pairs in the Cartesian product $X\times Y$. But to emphasize the point
about distinctions, we might employ the same trick of `marking' the elements
of $Y$, particularly when $Y=X$, with an asterisk. Then an alternative
construction of the product in $Sets$ is the set of \textit{unordered} pairs
$X\boxtimes Y=\left\{  \left\{  x,y^{\ast}\right\}  :x\in X;y^{\ast}\in
Y^{\ast}\right\}  $ which in the case of $Y=X$ would be $X\boxtimes X=\left\{
\left\{  x,x^{\ast}\right\}  :x\in X;x^{\ast}\in X^{\ast}\right\}  $. This
alternative construction of the product (isomorphic to the Cartesian product)
emphasizes the distinctions formed from $X$ and $Y$ so the ordering in the
ordered pairs of the usual construction $X\times Y$ is only a way to make the
same distinctions.

The set $X$ defines a partition $\pi_{X}$ on $X\times Y$ whose blocks are
$B_{x}=\left\{  \left(  x,y\right)  :y\in Y\right\}  =\left\{  x\right\}
\times Y$ for each $x\in X$, and $Y$ defines a partition $\pi_{Y}$ whose
blocks are $B_{y}=\left\{  \left(  x,y\right)  :x\in X\right\}  =X\times\{y\}$
for each $y\in Y$. Since $\pi_{X},\pi_{Y}\precsim\mathbf{1}_{X\times Y}$, the
induced maps (surjections if $X$ and $Y$ are non-empty) are the canonical
projections $p_{X}:X\times Y\rightarrow X$ and $p_{Y}:X\times Y\rightarrow Y$.

The universal mapping property for the product in $Sets$ is that given any
other `cone' of maps $f:Z\rightarrow X$ and $g:Z\rightarrow Y$, there is a
unique map $\left[  f,g\right]  :Z\rightarrow X\times Y$ such that the
triangles commute in the following diagram.

\begin{center}
$%
\begin{array}
[c]{ccccc}
&  & Z &  & \\
& \swarrow^{f} & ^{\exists!}\downarrow^{\left[  f,g\right]  } & ^{g}\searrow &
\\
X & \overset{p_{X}}{\longleftarrow} & X\times Y & \overset{p_{Y}%
}{\longrightarrow} & Y
\end{array}
$

Product diagram
\end{center}

From the data $f:Z\rightarrow X$ and $g:Z\rightarrow Y$, we need to
canonically construct the unique factor map $\left[  f,g\right]  :Z\rightarrow
X\times Y$. The map $f$ contributes the coimage $f^{-1}$ partition on $Z$ and
$g$ contributes the coimage $g^{-1}$ partition on $Z$ so we have the partition
lattice join $f^{-1}\vee g^{-1}$ in $\Pi\left(  Z\right)  $ whose blocks have
the form $f^{-1}\left(  x\right)  \cap g^{-1}\left(  y\right)  $. To define
the unique factor map $\left[  f,g\right]  :Z\rightarrow X\times Y$, the
discrete partition $\mathbf{1}_{Z}$ refines $f^{-1}\vee g^{-1}$ so there is a
canonical surjection $\mathbf{1}_{Z}\cong Z\twoheadrightarrow f^{-1}\vee
g^{-1}$which takes each $z\in Z$ to the unique block of the form
$f^{-1}\left(  x\right)  \cap g^{-1}\left(  y\right)  $ containing $z$. Each
such block of the form $f^{-1}\left(  x\right)  \cap g^{-1}\left(  y\right)  $
defines a mono $f^{-1}\vee g^{-1}\rightarrowtail X\times Y$ so the factor map
is: $\left[  f,g\right]  :Z\twoheadrightarrow f^{-1}\vee g^{-1}\rightarrowtail
X\times Y$.

\section{The Elements \& Distinctions Analysis of Equalizers and Coequalizers}

\subsection{The coequalizer in $Sets$}

For the equalizer and coequalizer, the data is not just two sets but two
parallel maps $f,g:X\rightrightarrows Y$. Then each element $x\in X$, gives us
a pair $f\left(  x\right)  $ and $g\left(  x\right)  $ so we take the
equivalence relation $\sim$ defined on $Y$ that is generated by $f\left(
x\right)  \sim g\left(  x\right)  $ for any $x\in X$. Then the coequalizer is
the quotient set $C=Y/\sim$ . When $\sim$ is represented as a partition on
$Y$, then it is refined by the discrete partition $\mathbf{1}_{Y}$ on $Y$, and
that refinement defines the canonical surjection $Y\cong\mathbf{1}%
_{Y}\rightarrow Y/\sim$.

For the UMP, let $h:Y\rightarrow Z$ be such that $hf=hg$. Then we need to show
there is a unique refinement-defined map $h^{\ast}:Y/\sim\rightarrow Z$ such
that the triangle commutes.

\begin{center}
$%
\begin{array}
[c]{ccccc}%
X & \underset{g}{\overset{f}{\rightrightarrows}} & Y & \twoheadrightarrow &
Y/\sim\\
&  &  & \searrow^{h} & ^{\exists!}\downarrow^{h^{\ast}}\\
&  &  &  & Z
\end{array}
$

Coequalizer diagram
\end{center}

We already have one partition $\sim$ on $Y$ which was generated by $f\left(
x\right)  \sim g\left(  x\right)  $. Since $hf=hg$, we are given that
$hf\left(  x\right)  =hg\left(  x\right)  $ so the coimage $h^{-1}$ has to
\textit{at least} identify $f\left(  x\right)  $ and $g\left(  x\right)  $
(and perhaps identify other elements) which means that $h^{-1}\precsim Y/\sim$
in the partition lattice on $Y$. Hence the induced surjection map $Y/\sim$
$\twoheadrightarrow h^{-1}$ and the mono $h^{-1}\cong h\left(  Y\right)
\rightarrowtail Z$ (taking $h^{-1}\left(  z\right)  $ to $z$) completes the
factor map $h^{\ast}:Y/\symbol{126}\twoheadrightarrow h^{-1}\rightarrowtail Z$.

\subsection{The equalizer in $Sets$}

For the same data $f,g:X\rightrightarrows Y$, the equalizer is the
$Eq=\left\{  x\in X:f\left(  x\right)  =g\left(  x\right)  \right\}  \subseteq
X$ so the map induced by that inclusion is the canonical map
$Eq\rightarrowtail X$.

The UMP is that for any other map $h:Z\rightarrow X$ such that $fh=gh,$then
$\exists!h_{\ast}:Z\rightarrow Eq$ such that the triangle commutes.

\begin{center}
$%
\begin{array}
[c]{ccccc}%
Eq & \rightarrowtail & X & \underset{g}{\overset{f}{\rightrightarrows}} & Y\\
^{\exists!}\uparrow^{h_{\ast}} & h\nearrow &  &  & \\
Z &  &  &  &
\end{array}
$

Equalizer diagram
\end{center}

\noindent The image of $h(Z)\subseteq X$ must satisfy $fh\left(  z\right)
=gh\left(  z\right)  $ for all $z\in Z$, so $f$ and $g$ agree on $h\left(
z\right)  \in X$ so $h\left(  Z\right)  \subseteq Eq$, which gives the
canonical injection $h\left(  Z\right)  \rightarrowtail Eq$ and the epi
$Z\twoheadrightarrow h\left(  Z\right)  $ completes the factor map $h_{\ast
}:Z\twoheadrightarrow h\left(  Z\right)  \rightarrowtail Eq$.

\section{The Elements \& Distinctions Analysis of Pushouts and Pullbacks}

\subsection{The pushout or co-Cartesian square in $Sets$}

It is a standard theorem of category theory that if a category has products
and equalizers, then it has all limits, and if it has coproducts and
coequalizers, then it has all colimits. Since we have presented the elements
\& distinctions analysis of the canonical maps for products and coproducts,
and for equalizers and coequalizers, the analysis extends to all limits and
colimits. Hence we have shown that the logical characterization of canonical
maps in $Sets$ satisfies Marquis's criteria:

\begin{quotation}
1. Morphisms that are part of the data of a limit are canonical morphisms; for
instance, the projection morphisms that are part of the notion of a product;

2. The unique morphism from a cone to a limit determined by a universal
property is a canonical morphism: and

3. In particular, the unique isomorphism that arise between two candidates for
a limit is a canonical morphism. \cite[p. 101]{marquis:canonical}
\end{quotation}

\noindent However, the theme would be better illustrated by considering some
more complicated limits and colimits such as Cartesian and co-Cartesian
squares, i.e., pullbacks and pushouts.

For the pushout or co-Cartesian square, the data are two maps $f:Z\rightarrow
X$ and $g:Z\rightarrow Y$ so we have the two parallel maps
$Z\overset{f}{\rightarrow}X\overset{i_{X}.}{\rightarrow}X\sqcup Y$ and
$Z\overset{g}{\rightarrow}Y\overset{i_{Y}.}{\rightarrow}X\sqcup Y$ and then we
can take their coequalizer $C$ formed by the equivalence relation $\sim$ on
the common codomain $X\sqcup Y$ which is the equivalence relation generated by
$x\sim y$ if there is a $z\in Z$ such that $f\left(  z\right)  =x$ and
$g\left(  z\right)  =y$. The canonical maps $X\rightarrow X\sqcup Y/\sim$ and
$Y\rightarrow X\sqcup Y/\sim$ are just the canonical injections into the
disjoint union followed by the canonical map of the coequalizer construction
analyzed above. As the composition of a canonical injection with a canonical
surjection, those canonical maps need not be either injective or surjective.

\begin{center}
$%
\begin{array}
[c]{ccccc}%
Z & \overset{f}{\rightarrow} & X & = & X\\
^{g}\downarrow^{{}} & \searrow & ^{{}}\downarrow^{can..} &  & \\
Y & \overset{can.}{\rightarrow} & C=X\sqcup Y/\sim &  & ^{{}}\downarrow^{h}\\
\shortparallel &  &  & ^{\exists!}\searrow^{h^{\ast}} & \\
Y &  & \overset{h^{\prime}}{\rightarrow} &  & U
\end{array}
$

Pushout or co-Cartesian square diagram
\end{center}

For the universal mapping property, consider any $h:X\rightarrow U$ and
$h^{\prime}:Y\rightarrow U$ such that $hf=h^{\prime}g$. Then $h^{-1}$ is a
partition on $X$ and $h^{\prime-1}$ is a partition on $Y$ so let $h^{-1}\sqcup
h^{\prime-1}$ be the disjoint union partition on $X\sqcup Y$. The condition
that for any $z\in Z$, $hf\left(  z\right)  =h^{\prime}g\left(  z\right)  =u$
for some $u\in U$ means that $h^{-1}\sqcup h^{\prime-1}$ must make at least
the identifications of the coequalizer (and perhaps more) so that
$h^{-1}\sqcup h^{\prime-1}$ is refined by $\sim$ as partitions on $X\sqcup Y$.
Since $h^{-1}\sqcup h^{\prime-1}\precsim\sim$ so each block $b$ in $\sim$ is
contained in a block of the form $h^{-1}\left(  u\right)  $ for some $u$ or a
block of the form $h^{\prime-1}\left(  u\right)  $ for some $u$. Hence that
block $b$ of $\sim$ is mapped by $h^{\ast}$ to the appropriate $u$ depending
on the case which defines the surjection from $X\sqcup Y/\sim$ to $h\left(
X\right)  \cup h^{\prime}\left(  Y\right)  \subseteq U$ and the inclusion
defines the injection to complete the definition of the canonical factor map
$h^{\ast}:C=X\sqcup Y/\sim\rightarrow U$.

\subsection{The pullback or Cartesian square in $Sets$}

For the Cartesian square or pullback, the data are two maps $f:X\rightarrow Z
$ and $g:Y\rightarrow Z$. We then have two parallel maps $X\times
Y\rightrightarrows Z$ (the projections followed by $f$ or $g$) so we take the
pullback as their equalizer $E$. The canonical maps $E\rightarrow X$ and
$E\rightarrow Y$ are the compositions of the canonical injective map
$E\rightarrow X\times Y$ followed by the canonical projections $p_{X}:X\times
Y\rightarrow X$ and $p_{Y}:X\times Y\rightarrow Y$. As the composition of a
canonical injection with a canonical surjection, those canonical maps need not
be either injective or surjective.

\begin{center}
$%
\begin{array}
[c]{ccccc}%
U &  & \overset{h}{\rightarrow} &  & X\\
& ^{\exists!}\searrow^{h_{\ast}} &  &  & \shortparallel\\
^{{}}\downarrow^{h^{\prime}} &  & E\subseteq X\times Y &
\overset{can.}{\rightarrow} & X\\
&  & \downarrow^{can.} &  & ^{{}}\downarrow^{f}\\
Y & = & Y & \overset{g}{\rightarrow} & Z
\end{array}
$

Pullback or Cartesian square diagram
\end{center}

For the universality property, consider any other maps $h:U\rightarrow X$ and
$h^{\prime}:U\rightarrow Y$ such that $fh=gh^{\prime}$. Hence $h^{\prime
}\left(  u\right)  $ and $h\left(  u\right)  $ are elements such that
$f\left(  h\left(  u\right)  \right)  =g\left(  h^{\prime}\left(  u\right)
\right)  $ so $\left(  h\left(  u\right)  ,h^{\prime}\left(  u\right)
\right)  \in E$ and thus for the images, there is the inclusion $h\left(
U\right)  \times h^{\prime}\left(  U\right)  \subseteq E$. Now $h$ contributes
the coimage partition $h^{-1}$ on $U$ and $h^{\prime}$ contributes the coimage
partition $h^{\prime-1}$ on $U$ and the join $h^{-1}\vee h^{\prime-1}$ is
refined by the discrete partition on $U$. Hence each $u\in U$ is contained in
a unique block $h^{-1}\left(  x\right)  \cap h^{\prime-1}\left(  y\right)  $
of the join so the refinement-induced canonical map $U\twoheadrightarrow
h\left(  U\right)  \times h^{\prime}\left(  U\right)  \subseteq E$ is defined
by $u\longmapsto\left(  x,y\right)  $ and the inclusion-defined injection
$h\left(  U\right)  \times h^{\prime}\left(  U\right)  \rightarrowtail E$
completes the definition of the canonical factor map $h_{\ast}:U\rightarrow E$.

\section{The Elements \& Distinctions Analysis of the Exponential Adjunction}

The adjunction $Hom_{Sets}\left(  X\times Y,Z\right)  \cong Hom_{Sets}\left(
X,Hom\left(  Y,Z\right)  \right)  $ is entirely in $Sets$. In the following diagram:

\begin{center}
$%
\begin{array}
[c]{ccccc}%
X & \overset{can.}{\rightarrow} & Hom\left(  Y,X\times Y\right)  &  & X\times
Y\\
& \searrow^{g} & ^{Hom\left(  Y,f\right)  }\downarrow & \exists!f &
f\downarrow\\
&  & Hom\left(  Y,Z\right)  &  & Z
\end{array}
$
\end{center}

\noindent the two canonical maps that need to be analyzed are unit
$X\overset{can.}{\rightarrow}Hom\left(  Y,X\times Y\right)  $ and the unique
factor map $f:X\times Y\rightarrow Z$ where the given data for the factor map
is the map $g:X\rightarrow Hom\left(  Y,Z\right)  $. That given data $g$
defines the partition $B_{z}=\left\{  \left(  x,y\right)  \in X\times
Y:g\left(  x\right)  \left(  y\right)  =z\right\}  $ on $X\times Y$ indexed by
$z\in Z$ and then the factor map $f:X\times Y\rightarrow Z$ is induced by the
refinement $\left\{  B_{z}\right\}  _{z\in Z}\precsim1_{X\times Y}$.

To analyze $X\overset{can.}{\rightarrow}Hom\left(  Y,X\times Y\right)  $, for
each $x\in X$, there is a (discrete) partition on $Y$ with the block $\left\{
y\right\}  $ indexed by the ordered pair $\left(  x,y\right)  $. Since the
discrete partition on $Y$ refines itself, there is the induced map
$h_{x}:Y\rightarrow X\times Y$ defined by $y\longmapsto\left(  x,y\right)  $.
And those functions $h_{x}\in Hom\left(  Y,X\times Y\right)  $ can also index
the blocks $\left\{  x\right\}  $ of the discrete partition on $X$ which
induces the injection $X\rightarrow Hom\left(  Y,X\times Y\right)  $ defined
by $x\longmapsto h_{x}$.

In the other UMP diagram for the adjunction:

\begin{center}
$%
\begin{array}
[c]{ccccc}%
Z & \overset{\operatorname{eval}}{\longleftarrow} & Hom\left(  Y,Z\right)
\times Y &  & Hom\left(  Y,Z\right) \\
& \nwarrow^{f} & g\times1_{Y}\uparrow & \exists!g & ^{g}\uparrow\\
&  & X\times Y &  & X
\end{array}
$
\end{center}

\noindent the two canonical maps that need to be analyzed are counit
$Z\overset{\operatorname{eval}}{\longleftarrow}Hom\left(  Y,Z\right)  \times
Y$ and the unique factor map $g:X\rightarrow Hom\left(  Y,Z\right)  $ where
the given data is the map $f:X\times Y\rightarrow Z$. The given data $f$
defines for each $x\in X$, a $Y$-partition $B_{z}=\left\{  y:f\left(
x,y\right)  =z\right\}  $ indexed by $z\in f\left(  X\times Y\right)  $ which
is refined by the discrete partition on $Y$ . Thus each $x\in X$ determines a
function $g_{x}:Y\rightarrow Z$, so the factor map is $g:X\rightarrow
Hom\left(  Y,Z\right)  $ where $g\left(  x\right)  =g_{x}\in Hom\left(
Y,Z\right)  $.

To analyze the canonical evaluation map $Hom\left(  Y,Z\right)  \times
Y\rightarrow Z$, each $z\in Z$ determines a partition on the domain by the
blocks $B_{z}=\left\{  \left(  h,y\right)  :h\left(  y\right)  =z\right\}  $,
and that partition is refined by the discrete partition on $Hom\left(
Y,Z\right)  \times Y$ and the induced map is the evaluation map.

\section{Example: A more complex canonical map}

Marquis \cite{marquis:canonical} gives the standard examples of canonical maps
that arise from limits and colimits but also mentions a more complex example
that will be analyzed. Let $\mathcal{C}$ be a category with finite products,
finite coproducts, and a null object (an object that is both initial and
terminal). Then a canonical morphism can be constructed from the coproduct of
two (or any finite number of) objects to the products of the objects: $X\sqcup
Y\rightarrow X\times Y$. In such a category abstractly specified, the map
could be constructed from the `atomic' canonical morphisms that are already
given by the arrow-theoretic definitions of products, coproducts, and the null
object. But the its \& dits analysis shows how all these `atomic' canonical
morphisms and their `molecular' compositions are not just assumed but are
constructed in $Sets$ or $Sets$-based categories according to the logical
theory of canonical maps.

There is a simple $Sets$-based category that has finite products, finite
products, and a null object, namely the category $Sets_{\ast}$ of pointed sets
where the objects are sets with a designated element (or basepoint), e.g.,
$\left(  X,x_{0}\right)  $ with $x_{0}\in X$, and the morphisms are set
functions that preserve the basepoints. The designation of the basepoint can
be given by a set map $1\overset{x_{0}}{\rightarrow}X$ in $Set$ which is taken
as part of the structure and is thus assumed canonical in $Set_{\ast}.$ A
basepoint preserving map $\left(  X,x_{0}\right)  \rightarrow\left(
Y,y_{0}\right)  $ is a set map $X\rightarrow Y$ in $Sets$ so that the
following diagram commutes:

\begin{center}
$%
\begin{array}
[c]{ccc}%
1 &  & \\
^{{}}\downarrow^{x_{0}} & \searrow^{y_{0}} & \\
X & \rightarrow & Y
\end{array}
$.
\end{center}

\noindent Hence $Sets_{\ast}$ can also be seen as the slice category $1/Sets$
of $Sets$ under $1$.

The null object is `the' one-point set $1$ and instead of assuming the
canonical morphisms that make it both initial and terminal, we need to
construct them using the its \& dits analysis. We have already seen that the
refinement relation $\mathbf{0}_{X}\precsim\mathbf{1}_{X}$ induces the unique
map $X\rightarrow1$ that makes $1$ the terminal object in $Sets$. And since
$1\overset{x_{0}}{\rightarrow}X\rightarrow1=1\overset{id.}{\rightarrow}1$, it
is also the terminal object in $Sets_{\ast}$. Moreover, the basepoint in
$(Y,y_{0})$ is given by the structurally canonical map $1\overset{y_{0}%
}{\rightarrow}Y$ and since $1\overset{id.}{\rightarrow}1\overset{y_{0}%
}{\rightarrow}Y=1\overset{y_{0}}{\rightarrow}Y$, that map $1\overset{y_{0}%
}{\rightarrow}Y$ is the unique map that makes $1$ also the initial object in
$Sets_{\ast}$. Hence in $Sets_{\ast}$, there is always a canonical map formed
by the composition: $X\rightarrow1\rightarrow Y=X\rightarrow Y$ (called the
\textit{zero arrow}).

To build up the its \& dits analysis of the canonical morphism $X\sqcup_{\ast
}Y\rightarrow X\times_{\ast}Y$ from the coproduct to the product in
$Sets_{\ast}$, we begin with the construction of the coproduct $X\sqcup_{\ast
}Y$ which is just the pushout in $Sets$ of the two canonical basepoint maps:

\begin{center}
$%
\begin{array}
[c]{ccccc}%
1 & \overset{x_{0}}{\rightarrow} & X & = & X\\
^{y_{0}}\downarrow^{{}} & \searrow & ^{{}}\downarrow^{can..} &  & \\
Y & \overset{can.}{\rightarrow} & X\sqcup_{\ast}Y=X\sqcup Y/\sim &  & ^{{}%
}\downarrow^{h}\\
\shortparallel &  &  & ^{\exists!}\searrow^{h_{\ast}} & \\
Y &  & \overset{h^{\prime}}{\rightarrow} &  & U
\end{array}
$
\end{center}

Since the only points in $X$ and $Y$ that are the image of an elements in $1$
are the basepoints, the equivalence relation $\sim$ only identifies the
basepoints $x_{0}$ and $y_{0}$. Hence $X\sqcup_{\ast}Y$ is like $X\sqcup Y$
except that the two basepoints are identified in the quotient $X\sqcup_{\ast
}Y=X\sqcup Y/\sim$ and that block identifying the basepoints is the basepoint
of $X\sqcup_{\ast}Y$. Then for any two set maps $h:X\rightarrow U$ and
$h^{\prime}:Y\rightarrow U$ that are also $Sets_{\ast}$ morphisms (i.e.,
preserve basepoints), there is a unique canonical factor map $h_{\ast}%
:X\sqcup_{\ast}Y\rightarrow U$ by the UMP for the pushout in $Sets$ to make
the triangles commute (and thus preserve basepoints). Hence $X\sqcup_{\ast}Y$
is the coproduct in $Sets_{\ast}$.

In a similar manner, one shows that the product $X\times_{\ast}Y$ in
$Sets_{\ast}$ is just the product $X\times Y$ in $Sets$ with $\left\langle
x_{0},y_{0}\right\rangle $ as the basepoint. Using the UMP of the product
$X\times_{\ast}Y$ in $Sets_{\ast}$, the two $Sets_{\ast}$ maps $1_{X}%
:X\rightarrow X$ and the canonical $X\rightarrow1\rightarrow Y$, we have the
unique canonical factor map $X\rightarrow X\times_{\ast}Y$ in $Sets_{\ast}$
and similarly for $Y\rightarrow X\times_{\ast}Y$ in $Sets_{\ast}$.

Then we put all the canonical maps together and use the UMP for the coproduct
in $Sets_{\ast}$ to construct the desired canonical map: $X\sqcup_{\ast
}Y\rightarrow X\times_{\ast}Y$ in $Sets_{\ast}$.

\begin{center}
$%
\begin{array}
[c]{ccccc}%
X & \overset{can.}{\longrightarrow} & X\sqcup_{\ast}Y &
\overset{can.}{\longleftarrow} & Y\\
& \searrow^{can.} & ^{\exists!}\downarrow^{can.} & ^{can.}\swarrow & \\
&  & X\times_{\ast}Y &  &
\end{array}
$

Coproduct diagram in $Sets_{\ast}$
\end{center}

This example shows how in a $Sets$-based category like $Sets_{\ast}$, the
given canonical maps for the structured sets (i.e., the basepoint maps
$1\overset{x_{0}}{\rightarrow}X$) are combined with the canonical maps defined
by the its \& dits analysis in $Sets$ to give the canonical morphisms in the
$Sets$-based category. In abstract category theory, as in the case of a
category $\mathcal{C}$ which is assumed to have finite products, finite
coproducts, and a null object, the `atomic' canonical morphisms are all given
as part of the assumed UMPs for products, coproducts, and the null object
which are then composed to define other `molecular' canonical morphisms. That
suggests that the only categories where a \textit{theory} of canonical
morphisms is needed is $Sets$ and $Sets$-based categories, and that is the
theory presented here.

\section{Concluding philosophical reflections}

The "logical" in the logical theory of canonicity refers to the two dual
mathematical logics: the Boolean logic of subsets and the logic of partitions.
Note that from the mathematical viewpoint, the Boolean logic of subsets and
the logic of partitions have \textit{equal} intertwining roles in the whole
analysis. Normally, we might say that "subsets" and "partitions" are
category-theoretic duals, but we have tried to show a more fundamental
analysis based on "elements \& distinctions" or "its \& dits" that are the
building blocks of subsets and partitions and that \textit{underlie} the
duality in $Sets$.

Thus instead of saying that duality explains elements and distinctions, we
tried to show that the intricate and precise interplay of elements and
distinctions explains morphisms, duality, canonicity, and universal
constructions in $Sets$, which generalizes to other $Sets$-based concrete
categories and which is abstracted in abstract category theory.

Our focus here is the E\&D treatment of canonicity.

\begin{itemize}
\item Each construction starts with certain data.

\item When that data is sufficient to define inclusions in an associated
subset lattice or refinements in an associated partition lattice, then the
induced `logical' maps (and their compositions) are canonical.
\end{itemize}

This suggests that the dual notions of elements \& distinctions (its \& dits)
have some broader significance. One possibility is they are respectively
mathematical building blocks of the old metaphysical concepts of matter (or
substance) and form (as in in-form-ation). The matter versus form idea
\cite{aimsworth:FandM} can be illustrated by comparing the two lattices of
subsets and partitions on a set--the two lattices that we saw defined the
canonical morphisms and canonical factor maps in $Sets$-based categories.

For $U=\left\{  a,b,c\right\}  $, start at the bottom and move towards the top
of each lattice.%

\begin{center}
\includegraphics[
height=1.5826in,
width=5.6896in
]%
{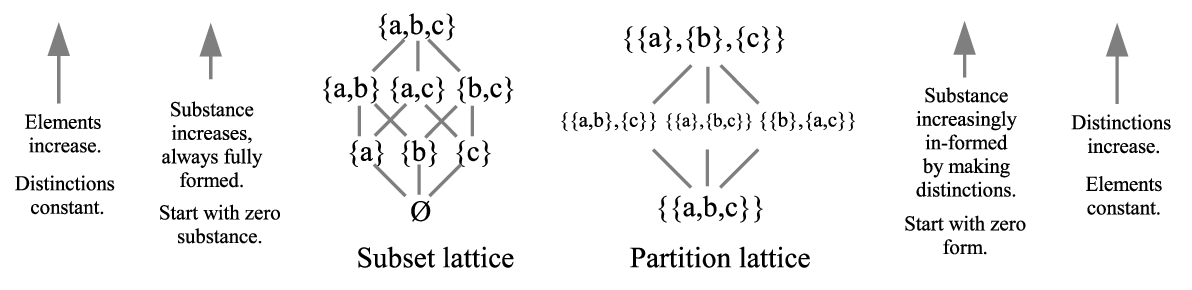}%
\end{center}

\begin{center}
Figure 3: Moving up the subset and partition lattices.
\end{center}

At the bottom of the Boolean subset lattice is the empty set $\emptyset$ which
represents no substance (no elements or `its'). As one moves up the lattice,
new elements of substance, new elements, are created that are always fully
distinguished or formed until finally one reaches the top, the universe $U$.
Thus new substance is created in moving up the lattice but each element is
fully formed and thus distinguished from the other elements.

At the bottom of the partition lattice is the indiscrete partition or "blob"
$\mathbf{0}_{U}=\left\{  U\right\}  $ (where the universe set $U$ makes one
block) which represents all the substance or matter but with no distinctions
to in-form the substance (no distinctions or `dits'). As one moves up the
lattice, no new substance is created but distinctions are created that in-form
the indistinct elements as they become more and more distinct. Finally one
reaches the top, the discrete partition $\mathbf{1}_{U}$, where all the
elements of $U$ have been fully in-formed or distinguished. A partition
combines indefiniteness (within blocks) and definiteness (between blocks). At
the top of the partition lattice, the discrete partition $\mathbf{1}%
_{U}=\left\{  \left\{  u\right\}  :\left\{  u\right\}  \subseteq U\right\}  $
is the result making all the distinctions to eliminate any indefiniteness.
Thus one ends up at essentially the same place (universe $U$ of fully formed
entities) either way, but by two totally different but dual `creation stories':

\begin{itemize}
\item Subset Creation Story: creating elements as in creating fully-formed and
distinguished matter out of nothing, versus

\item Partition Creation Story: creating distinctions by starting with a
totally undifferentiated matter and then, in a `big bang,' start making
distinctions, e.g., breaking symmetries, to give form to the matter.
\end{itemize}

Moreover, we have seen that:

\begin{itemize}
\item the quantitative increase in substance (normalized number of elements)
moving up in the subset lattice is measured by logical or Laplacian
probability, and

\item the quantitative increase in form (normalized number of distinctions)
moving up in the partition lattice is measured by logical entropy
(\cite{ell:countingdits}; \cite{ell:nf4it}; \cite{manf:intro}).
\end{itemize}

Declarations:

The author has no competing or conflicting interests to declare that are
relevant to the content of this article.

No funds, grants, or other support was received.

Data availability: N/A

\end{document}